\documentstyle[12pt]{article}
\newtheorem{th}{Theorem}[section]
\newtheorem{pr}[th]{Proposition}
\newtheorem{lem}[th]{Lemma}
\newtheorem{co}[th]{Corollary}
\newtheorem{remark}[th]{Remark}
\newtheorem{example}[th]{Example}
\newtheorem{definition}[th]{Definition}
\newenvironment{proof}{{\noindent {\em Proof:}}}{\hspace*{1cm}
        $\Box$\medskip}
\newenvironment{re}{\begin{remark}\rm}{\end{remark}}

\newenvironment{de}{\begin{definition}\rm}{\end{definition}}

\def\2stack#1#2{\mathrel{\mathop{#1}\limits_{#2}}}
\def\3stack#1#2#3{\mathrel{\mathop{\mathop{#1}\limits_{#2}}\limits_{#3}}}
\def\vbar{\mathchoice{\vrule height6.3ptdepth-.5ptwidth.8pt\kern-.8pt}
   {\vrule height6.3ptdepth-.5ptwidth.8pt\kern-.8pt}
   {\vrule height4.1ptdepth-.35ptwidth.6pt\kern-.6pt}
   {\vrule height3.1ptdepth-.25ptwidth.5pt\kern-.5pt}}
\def\fudge{\mathchoice{}{}{\mkern.5mu}{\mkern.8mu}}
\def\bbc#1#2{{\rm \mkern#2mu\vbar\mkern-#2mu#1}}
\def\bbb#1{{\rm I\mkern-3.5mu #1}}
\def\bba#1#2{{\rm #1\mkern-#2mu\fudge #1}}
\def\bb#1{{\count4=`#1 \advance\count4by-64 \ifcase\count4\or\bba
A{11.5}\or
   \bbb B\or\bbc C{5}\or\bbb D\or\bbb E\or\bbb F \or\bbc G{5}\or\bbb H\or
   \bbb I\or\bbc J{3}\or\bbb K\or\bbb L \or\bbb M\or\bbb N\or\bbc O{5} \or
   \bbb P\or\bbc Q{5}\or\bbb R\or\bbc S{4.2}\or\bba T{10.5}\or\bbc U{5}\or
   \bba V{12}\or\bba W{16.5}\or\bba X{11}\or\bba Y{11.7}\or\bba Z{7.5}\fi}}

\newcommand{\C}{{\bb C}}
\newcommand{\PP}{{\bb P}}
\newcommand{\r}[1]{~\mbox{$(${\rm \ref{#1}}$)$}}

\begin{document}
\title{DEGREE OF THE GENERALIZED PL\"UCKER EMBEDDING
OF A QUOT SCHEME AND QUANTUM COHOMOLOGY}
  \author{M.S. Ravi
\thanks{East Carolina University, Greenville, NC 27858.}
 \and J. Rosenthal
\thanks{University of Notre Dame, Notre Dame, IN 46556.
This author was supported  in part by
 NSF grant DMS-9201263.}
\and X. Wang
\thanks{Texas Tech University, Lubbock, TX 79409.
This author was supported in part by NSF grant
DMS-9224541.}  }
\date{}\maketitle

\begin{abstract} We compute the degree of the generalized
Pl\"ucker embedding $\kappa$ of a Quot scheme $X$  over $\PP^1$. The
space $X$ can also be considered as a compactification of the space
of algebraic maps of  a fixed degree from $\PP^1$ to the Grassmanian
$\rm{Grass}(m,n)$. Then the degree of the embedded variety
$\kappa (X)$ can be interpreted as an intersection
product of pullbacks of cohomology classes from
$\rm{Grass}(m,n)$ through the map $\psi$ that evaluates a
map from $\PP^1$ at a point $x\in \PP^1$.  We show that
our formula for the degree verifies the formula
for these intersection products predicted by physicists
through Quantum cohomology~\cite{va92}~\cite{in91}~\cite{wi94}.
We arrive at the degree by
proving a version of the classical Pieri's formula on the variety $X$,
using a cell decomposition of a space that lies in between
$X$ and $\kappa (X)$.
\end{abstract}

Let $X'$ be the space of all algebraic maps, of a fixed degree $q$,
from the
projective line $\PP^1$ to $\mbox{Grass}(m,n)$,
the Grassmannian of all
$m$-dimensional subspaces of a fixed $n$-dimensional vector
space  $V$. In various papers, physicists have been discussing so called
correlation functions on this
space~(\cite{in91}~\cite{va92}~\cite{wi94}) which specify the
intersection products on $X'$ of  pullbacks of cohomology classes from the
Grassmannian (a precise formulation can be found in Section 4). Further, based
on certain physical arguments, there  have been some conjectured formulas for
these correlation functions (actually these conjectures deal with the more
general case of maps from  any Riemann surface to $\mbox{Grass}(m,n)$, but we
will restrict our attention to maps from $\PP^1$).
In~\cite{be93} the authors give
a mathematically rigorous proof of the conjectured formula in the
case of maps from a Riemann surface of genus one to $\mbox{Grass}(2,n)$.
In this paper, we show that the conjecture holds for a certain
class of intersection products on the space of maps
from $\PP^1$ into $\mbox{Grass}(m,n)$ for any $m$ and $n$~\r{phythm}.
Some intriguing connections with quantum cohomology and Floer cohomology
can be found in~\cite{gi94}~\cite{sad}~\cite{as}~\cite{pi94}.

We shall be working on a natural compactification of $X'$,
namely a Quot scheme defined as follows: Let
$V_{\PP}=V\otimes O_{\PP^1}$ and
let $f$ be the polynomial $f(l)= p(l+1)+q$, where $p=n-m$.
As explained in~\cite{st87},
the quotient scheme $X=\mbox{Quot}^f_{V_\PP}$ that
parameterizes all quotient sheaves $\cal B$ of $V_{\PP}$ such that the
Hilbert polynomial $\chi ({\cal B}(l))=f(l)$, can be considered as
a compactification of $X'$. Some technical points on the choice of this
compactification and the independence of the physical predictions
on this choice can be
found in Section 4. The space $X$ can be  concretely described as
follows~(\cite{ra94}): each point in $X$ can be considered as an equivalence
class of matrices  $
M=(M_{ij}(s,t))_{1\le i\le m,\ 1\le j\le n}
$
where each $M_{ij}(s,t)$ is a homogeneous polynomial of degree $q_i$
and $\sum q_i=q$. Two such matrices $M$ and $M'$ are considered
equivalent if, after rearranging the rows if necessary, the row degrees
of $M$ and $M'$ are the same and there exists an $m\times m$
matrix $\, U$ such that its entry $U_{ij}(s,t)$ is a homogeneous
polynomial
of degree $q_j-q_i$, $\det U$ is a non-zero constant and $M'=UM$.
In the sequel we will find it convenient to switch back and forth
between this and the description of the points in $X$ as quotient
sheaves.

Given a sheaf ${\cal B}\in X$, let
$$
0\to {\cal A}\to V_{\PP}\to {\cal B}\to 0
$$
be the exact sequence defining $\cal B$. If we choose
an isomorphism ${\cal A}\simeq
{m\atop{\oplus\atop{i=1}}}
O_{\PP}(-q_i)$ and choose basis for $V$, then the
quotient sheaf $\cal B$ can be identified with the
matrix $M$ of the map ${\cal A}\to V_{\PP}$.

There is a  natural map from $X$ into a projective space,
namely the generalized Pl\"ucker map $\kappa$ defined as follows:
 to each matrix $M\in X$, $\kappa$ assigns the $m\times m$ minors of
$M$, considered as a point in
$$
\PP=\PP(\wedge^m V\otimes H^0(\PP^1,O_{\PP^1}(q))).
$$
The map $\kappa$ is not an embedding and the image, which we denote by
$K^q_{m,p}$,  is a singular variety.
In this paper we compute the degree of $K^q_{m,p}$ and that
of certain
subvarieties of this space that can be thought of as pullbacks of Schubert
varieties through the evaluation map from $X$ to
$\mbox{Grass}(m,n)$ (see Section 4 for a precise formulation). These
degrees are a particular case of the conjectures from physics and we show that
our formulas agree with the predicted degrees.
Since the degree of $K^q_{m,p}$ is also equal to the degree
of the pole placement map in the critical dimension we derive
in this paper also an important result in Systems theory.

Our methods are quite elementary.
We generalize the classical methods used to
compute the degree of the Grassmannian. We use an intermediate space,
$A^q_{m,p}$
such that the map $\kappa$ from $X$ to $\PP$ factors as
\begin{equation}\label{plu}
X\stackrel{\phi}{\rightarrow}A^q_{m,p}\stackrel{\pi}{\rightarrow}\PP.
\end{equation}
The map $\phi$
sends a matrix $M(s,t)\in X$ to $M(s,1)$ and the map $\pi$ sends
$M(s,1)$ to its $m\times m$ minors homogenized to polynomials of
total degree $q$.  The key ingredient for our work is
a cellular decomposition of $A^q_{m,p}$ described in~\cite{wa94}
that parallels the standard cellular decomposition of a Grassmannian,
though this situation is more complicated. This cellular decomposition
enables us to use B\'{e}zout's theorem on $\PP$ to prove an
analogue of Pieri's formula (Proposition~\ref{zint}) on $K^q_{m,p}$.
This reduces the
computation of the degree to a combinatorial problem.

In Section~1 we define the space $A^q_{m,p}$ and recall its basic
properties and its cell decomposition (Proposition~\ref{c}).
In Section~2 we define the
subvarieties of $K^q_{m,p}$ (Definition~\ref{zdef})
that play a role analogous to that
of the Schubert subvarieties of a Grassmannian. In  Section 3
 we apply
B\'{e}zout's theorem to the intersection of these subvarieties by appropriate
hyperplanes of $\PP$. We then give
a combinatorial description for the degree of $K^q_{m,p}$
and its subvarieties.
In Section 4 we give a precise formulation of the Physics conjecture and
interpret the degrees we compute in terms of this conjecture.
We then show that
our computation of the degrees  agrees with
the  formula conjectured by the physicists.

We started this paper with a view of understanding the space
$X$, considered as a compactification of the space of
all $m$-input, $p$-output transfer functions of McMillan degree
$q$~(\cite{ra94},\cite{ro94}) in Systems theory. A substantial portion
of this paper has been in circulation
for some time under the title ``Degree of the Generalized Pl\"ucker
embedding of a
Quot scheme".

\section{The cellular decomposition of $A^q_{m,p}$}
\setcounter{equation}{0}

In this section we summarize those
 results obtained in~\cite{wa94} which we will use
in this paper. The space $A^q_{m,p}$ consists of equivalence
classes of polynomial matrices
$M(s)=(M_{ij}(s))_{1\le i\le m,\ 1\le j\le n}$
such that the degree of any $m\times m$
minor of $M$ is at most $q$ and at least one of these minors is
a non-zero polynomial. Two such
matrices $M(s)$ and $M'(s)$ are considered equivalent, if there exists a
unimodular polynomial matrix $U(s)$  such that $M'(s)=U(s)M(s)$.

\begin{de}
Given any $m\times n$ polynomial matrix
$M(s)$, there exist
unique $\nu=(\nu_{1},\dots,\nu_{m})$ with $\nu_{1}\leq \cdots  \leq
\nu_{m}$ and
$$\sum_{i=1}^{m}\nu_{i}=
\mbox{maximum degree of $m\times m$ minors of $M(s)$} $$
and an $m\times m$
unimodular matrix $U(s)$ such that the matrix $M'(s)=U(s)M(s)$
has row degrees $\nu_{1}\leq \cdots  \leq
\nu_{m}$.
The numbers $\nu =(\nu_{1}, \ldots , \nu_{m})$ are
called the ordered Kronecker indices of the equivalence class of
$M(s)$.
\end{de}

A matrix $M(s)$ is called row reduced if
the ordered Kronecker indices are equal to the degrees of
the rows of $M(s)$.
$M(s)$ is row reduced if and only if the high order coefficient matrix
of $M(s)$ has full rank, where the high order coefficient matrix
of a polynomial matrix $M(s)$ is a matrix whose entries of the $i$-th row
are the coefficients of $s^{\nu_i}$ of the $i$-th row of $M(s)$ where
$\nu_i$ is the highest power of $s$ in the $i$-th row of $M(s)$.

\begin{de}(\cite{fo75}\cite{wa94})
Given an $m\times n$ row reduced polynomial matrix $M(s)$ with $M_h$
the high order coefficient matrix of $M(s)$,
the $i$-th pivot index $\mu_i'$ is the largest integer such
that the submatrix of $M_h$ formed from the intersection of columns
$\mu_1',\dots,\mu_i'$ with the rows corresponding to
indices $\leq \nu_i$ has rank $i$. The ordered  pivot indices
$\mu=(\mu_1,\mu_2,\dots,\mu_m)$ of the equivalence class of
$M(s)$
are the indices obtained from
$(\mu_1',\mu_2',\dots,\mu_m')$ by reordering such that $\mu_i<\mu_{i+1}$
if $\nu_i=\nu_{i+1}$.
\end{de}
Let $\tilde{I}(m,p)$ be the set of $m$-tuple of integers defined by
$$
\tilde{I}(m,p)=\{\alpha=(\alpha_1,\dots,\alpha_m)|
1\le \alpha_1<\cdots<\alpha_m,\ \alpha_k\neq \alpha_l \bmod n\mbox{ for }
k\ne l\}.
$$

\begin{de} For each equivalence class $M\in A^q_{m,p}$ with Kronecker
index $\nu =(\nu_1,\dots ,\nu_m)$ and ordered pivot index
$\mu=(\mu_1,\dots ,\mu_m)$ we define a new index
$\alpha=(\alpha_1,\dots \alpha_m)\in \tilde{I}(m,p)$
by
$\alpha_l :=n\nu_l+\mu_l$.
\end{de}
Further for each index $\alpha\in \tilde{I}(m,p)$ we define
\begin{equation}\label{norm2}
|\alpha|=\sum_{l=1}^m (\alpha_l-l)-
\sum_{l=2}^m\sum_{k=1}^{l-1}\left[ {\frac{\alpha_l-\alpha_k}{n}}\right]
\end{equation}
where $[r]$ is the largest integer less than or equal to $r$.
We also define a partial order on $\tilde{I}(m,p)$ as follows:
first associate with each
$\alpha=(\alpha_1,\dots,\alpha_m)$ an infinite
sequence:
$$
f(\alpha)=(f_1(\alpha),f_2(\alpha),\cdots)
$$
where
$$
\{ f_l(\alpha) \}
=\{\alpha_j+k(n)\mid j=1,\dots,m\mbox{ and }k=0,1,2,\dots,\}
$$
and  order the set such that
$$f_1(\alpha)<f_2(\alpha)<\cdots,$$
and then define the partial order on $\tilde{I}(m,p)$ by
\begin{equation}\label{partial2}
\alpha\leq \beta\mbox{ if and only if }
f_l(\alpha)\leq f_l(\beta)\mbox{ for all $l$}.
\end{equation}

Next we define a topology on $A^q_{m,p}$. Let ${\cal P}^q_{m,p}$
be the set of  all $m\times n$ full rank polynomial
matrices of degree at most $q$ and denote with
${\cal P}^{q,r}_{m,p}$ the subset of ${\cal P}^q_{m,p}$  formed by all
matrices whose entries are polynomials of degree  at most $r$.
Then
$${\cal P}^{q,0}_{m,p}\subset {\cal P}^{q,1}_{m,p}\subset
{\cal P}^{q,2}_{m,p}\subset \cdots$$
with the union $${\cal P}^q_{m,p}=
\bigcup_{r=0}^{\infty}{\cal P}^{q,r}_{m,p}.$$
The set of all $m\times n$ polynomial matrices whose entries
are polynomials of degree at most $r$ is an affine space
$$\C^{mn(r+1)}$$
and the conditions that the degrees of the $m\times m$ minors are
at most $q$ are polynomial conditions on $\C^{mn(r+1)}$
which defines an algebraic set.
${\cal P}^{q,r}_{m,p}$ is a Zariski open set of this algebraic set.
Take the Zariski topology on ${\cal P}^{q,r}_{m,p}$.
The direct limit of the topologies on
${\cal P}^{q,r}_{m,p}\ \ r=0,1,\dots$ defines a topology
on ${\cal P}^q_{m,p}$. In other
words, a subset of ${\cal P}^q_{m,p}$ is open  if and only if its
intersection with ${\cal P}^{q,r}_{m,p}$ is open as a subset of
${\cal P}^{q,r}_{m,p}$ for each $r$. The topology which
we take on $A^q_{m,p}$ is the quotient
topology under  row  equivalence,
i.e.
a subset $U$ of $A^q_{m,p}$ is open if, and only if
the subset $V$ of $P^q_{m,p}$ formed by all the polynomial
matrices in the equivalence classes of $U$ is open.
Since the minors of an $m\times n$ matrix are polynomials of its
entries,
any polynomial condition on $K^q_{m,p}$ induces a polynomial
condition on ${\cal P}^q_{m,p}$. So the map $\pi$ defined in\r{plu} is
continuous under the topology we defined on $A^q_{m,p}$.

Let $C_\alpha$ be the subset of $A^q_{m,p}$ consisting all the
elements with index $\alpha$. The main result given
in~\cite{wa94} can then be summarized in the following
proposition:
\begin{pr}~\cite{wa94} \label{c}
\begin{itemize}
\item[1)] $C_\alpha$ is an open cell of dimension $|\alpha|$.
\item[2)] $C_\alpha\cap C_\beta=\emptyset$ if $\alpha\neq
\beta$.
\item[3)] $\overline{C}_\alpha={\displaystyle
\bigcup_{\beta\in \tilde{I}(m,p) \atop \beta\leq\alpha}}C_\beta$.
\end{itemize}
\end{pr}

\section{Generalized Schubert Subvarieties of $K^q_{m,p}$}
\setcounter{equation}{0}

We fix the coordinates of $K^q_{m,p}$ first.
For each $i=(i_1,\dots,i_m)$, $0\le i_1<\cdots< i_m\le n$, let
$$
z_{(i;0)}t^q+z_{(i;1)}t^{q-1}s+\cdots+z_{(i;q)}s^q
$$
be the $m\times m$ minor of an $M(s,t)\in X$ consisting of the $i_1$th
through $i_m$th columns. Then
$$
z=(z_{(i;d)})_{0\le i_1<\cdots< i_m\le n,\ 0\le d \le q}.
$$
is the homogeneous coordinate of the image of $M$ in $K^q_{m,p}$.

Let $I(m)$ be the set of $m$-tuple of integers defined by
$$
I(m)=\{i=(i_1,\dots ,i_m)\,|\, 1\le i_1<\dots <i_m\}.
$$
Define
\begin{equation}\label{norm1}
|i|=\sum_{l=1}^m (i_l-l)
\end{equation}
and the partial order
\begin{equation} \label{partial1}
(i_1,\dots,i_m)\leq (j_1,\dots,j_m)\mbox{ if and only if }
i_l\leq j_l\mbox{ for all $l$}
\end{equation}
on $I(m)$.

Let $e_l$ be the vector whose $l$th component is $1$ and all the other
components are zero and
$$F_l={\rm sp}\ \{e_1,\dots,e_l\}.$$
For any $i=(i_1,\dots,i_m)\in I(m)$, $i\le (p+1,\dots,n)$, let $S_i$
be the Schubert variety $S(F_{i_1},\dots,F_{i_m})$ under the standard
flag $F_{i_1}\subset\cdots\subset F_{i_m}$; i.e.
\begin{equation}\label{schubert}
S_i=S(F_{i_1},\dots,F_{i_m})
=\{x\in{\rm Grass}(m,n)| \dim x\cap F_{i_l}\ge l\}.
\end{equation}
Then $S_i$ is a subvariety of dimension $|i|$ defined by
\begin{equation}\label{svariety}
S_i=\{x\in {\rm Grass}(m,n)
\mid x_j=0\mbox{ for all }j\not\leq i\}
\end{equation}
where $(x_i)$ are the Pl\"ucker coordinates of a point
$x\in{\rm Grass}(m,n)$.

We would like to generalize the Schubert varieties in
Grass$(m,n)=K^0_{m,p}$ to $K^q_{m,p}$.

For each $d\leq q$, the subset $\{z\in K^q_{m,p}\, |\, z_{(i:l)}=0\mbox{
for all }l>d\mbox{ and }i\}$ can be identified naturally with $K^d_{m,p}$.
By abuse of notation we will denote this subset as $K^d_{m,p}\subseteq
K^q_{m,p}$.

Let $\psi_d$ be the projection of
$z\in K^d_{m,p}\subseteq K^q_{m,p}$
on its components $(z_{(i;d)})$; i.e.
\begin{equation}\label{psidef}
\psi_d(z)=\{ z_{(i;d)}| 1\leq i_1<\cdots<i_m\leq n\}.
\end{equation}
Then $\psi_d$ is a rational map (it is in fact a central projection) and
$$
\psi_d(K^d_{m,p})={\rm Grass}(m,n).
$$

A first attempt to generalize Schubert varieties
could be to take $\psi_d^{-1}(S_i)$ in $K^d_{m,p}\subseteq K^q_{m,p}$
for each $d\le q$ and for each Schubert variety $S_i$
in ${\rm Grass}(m,n)$.
Unfortunately this is not a closed set.
A better definition comes from the cell decomposition of $A^q_{m,p}$.
But we first need the following definition.
\begin{de}
For each $(i;d)$, $i=(i_1,\dots,i_m)$, $1\leq i_1<\cdots<i_m\leq n$,
and $d=km+r$, $0\leq r<m$, let $\alpha =(i;d)$ where
$\alpha=(\alpha_1,\dots,\alpha_m)$
is defined through:
\begin{equation}\label{ialpha}
\alpha_l=\left\{\begin{array}{ll}
k(n)+i_{l+r}       & \mbox{for } l=1,2,\dots,m-r\\
(k+1)(n)+i_{l-m+r} & \mbox{for } l=m-r+1,\dots,m.
\end{array}\right.
\end{equation}
\end{de}

Then $C_\alpha$ with $\alpha=(\alpha_1,\dots,\alpha_m)=(i;d)$ is the
``thickest'' cell among all the cells $C_\beta$ in $\pi^{-1}(K^d_{m,p})$
such that $$\psi_d(\pi(C_\beta))\subset S_i.$$
It can be proved that $\pi(\overline{C}_\alpha)$ is a
subvariety (see Lemma~\ref{lem37} and Proposition~\ref{irreducible}) and
$\psi_d(\pi(\overline{C}_\alpha))=S_i$. So $\pi(\overline{C}_\alpha)$ is
the variety we want.
In order to give a definition similar to (\ref{svariety}) we
re-index the coordinates of ${\bb P}\supset K^q_{m,p}$,
as
$z_{(i;d)}=z_\alpha.$

We set
$$
 I(m,p)=\{\alpha=(\alpha_1,\dots,\alpha_m)\, |\, 1\le\alpha_1
<\cdots <\alpha_m,\ \alpha_m-\alpha_1<n\}.
$$
If  $i\in I(m)$ and  $i\leq (p+1,\dots,n)$ then $\alpha=(i;d)\in I(m,p)$
and for each $\alpha\in I(m,p)$ there exists  unique  $i\in I(m)$,
$i\le (p+1,\dots,n)$, and
$d$ such that $\alpha
= (i;d)$. Further, for $I(m,p)$
the notion of the partial order~(\ref{partial2}) agrees with the partial
order~(\ref{partial1}) and
the notion of $|\alpha|$ defined in~(\ref{norm2}) reduces
to~(\ref{norm1}). So $I(m,p)$ is a subset of $\tilde{I}(m,p)$ as well as
$I(m)$.

\begin{re}\label{remarki}
In a partially ordered set, an element $\alpha$ is said to
cover another element $\beta$ if $\beta<\alpha$ and there exists
no $\gamma$ such that $\beta<\gamma<\alpha$~\cite{kr86}.
{}From the definition of the
partial order one can see that for any
$\alpha=(i;d)$ and $\beta=(j;b)$ in $I(m,p)$,
$\alpha$ covers $\beta$ if and only if either
\begin{itemize}
\item[a.] $b=d$ and $i$ covers $j$ or
\item[b.] $b=d-1$, $i=(1,i_2,\dots,i_m)$ and $j=(i_2,\dots,i_m,n)$
for some $1<i_2<\cdots<i_m<n$.
\end{itemize}
\end{re}

\begin{de}\label{zdef}
For any $\alpha\in I(m,p)$, $\alpha=(i;d)\le (p+1,\dots n;q)$, let
$Z_\alpha$, or $Z_i^d$, be the closed subset of $K^q_{m,p}$
defined by $$
Z_\alpha=Z_i^d
=\{z\in K^q_{m,p}| z_\beta=0\mbox{ for all }\beta\not\leq\alpha \}
$$
and
$O_\alpha$, or $O_i^d$, be the open set of $Z_\alpha$ defined by
$$
O_\alpha=O_i^d
=\{z\in Z_\alpha| z_\alpha\neq 0 \}
$$
\end{de}

\begin{pr}\label{o}
$$
Z_\alpha=\bigcup_{\beta\in I(m,p), \atop \beta\leq\alpha} O_\beta.
$$
\end{pr}

\begin{proof}
Follows from the definition.
\end{proof}

\begin{lem}\label{lem37}
Let $\alpha=(i;d)$. Then
\begin{itemize}
\item[1)]
$
O_\alpha=\bigcup \pi(C_\beta),
$
where the union is over all cells in $A^q_{m,p}$
with  Kronecker indices
$
\nu=(\nu_1,\dots,\nu_m),\  \sum \nu_l=d
$ and  pivot index
$\mu=(\mu_1,\dots,\mu_m)$ such that
$\{\mu_1,\dots,\mu_m\}=\{i_1,\dots,i_m\}
$
as unordered sets.
\item[2)]
\begin{equation}\label{z}
Z_\alpha=\bigcup_{\beta\in \tilde{I}(m,p),
\atop \beta\leq \alpha}\pi(C_\beta)=\pi
(\overline{C}_\alpha). \end{equation}
\end{itemize}
\end{lem}

\begin{proof}
Using the echelon form of the elements in $C_\alpha$
(see~\cite{wa94}) one can see that
$
\pi (C_\alpha)\subset Z_\alpha.
$
Since
$C_\beta\subset \overline{C}_\alpha$
for all $\beta$ with
Kronecker indices
$\nu,\
\sum \nu_l =d$
and pivot indices as above
and since $\pi$ is continuous,
$
\bigcup \pi(C_\beta)\subset Z_\alpha.
$
Furthermore,
$z_\alpha\neq 0$
for all the points in $\pi(C_\beta)$.
So
$\bigcup \pi(C_\beta)\subset O_\alpha.
$

On the other hand, for any $z\in O_\alpha$ let $M(s)\in\pi^{-1} (z)$
be row reduced. Then by looking at the high order coefficient
matrix~\cite{fo75,wa94} of $M(s)$ one concludes immediately that
$M$ must have Kronecker index
$\nu,\
\sum \nu_l =d
$ and pivot index
$\mu$
such that $\{\mu_{1}, \ldots , \mu_{m}\}=\{i_1,\dots,i_m\}.
$

Therefore
$\bigcup \pi(C_\beta)=O_\alpha$
and (\ref{z}) follows because of Proposition~\ref{c} and
Proposition~\ref{o}.
\end{proof}

\begin{pr}\label{irreducible}
$Z_\alpha$ is an irreducible  subvariety of $K^q_{m,p}$ of dimension
$|\alpha|$.
\end{pr}

\begin{proof}
Since $\overline{C}_\alpha$ is
irreducible and $\pi$ is one to one on the open set consisting
 of all matrices whose maximal minors do not have any common factors,
$Z_\alpha=\pi(\overline{C}_\alpha)
$
is irreducible and
$\dim Z_\alpha=\dim
\overline{C}_\alpha=|\alpha|.
$ \end{proof}

\begin{re}
{}From Proposition~\ref{o} and Remark~\ref{remarki} one can see
immediately that
$$\psi_d (Z^d_i)=S_i$$
where $S_i$ is the Schubert variety of Grass$(m,n)$ defined
by~(\ref{svariety}). In the next section we shall show that
$Z_\alpha$ is birationally equivalent to the Schubert variety
$S_\alpha\subset {\rm Grass}(m,n(q+1))$
(see Proposition~\ref{pro41}). This is one of the
reasons why we use two kinds of indices to label the variety.
When $d=0$, $Z^0_i$ reduces to the Schubert variety~$S_i$.
\end{re}

\section{Degree of $Z_\alpha$}
\setcounter{equation}{0}

To prove our formula for the degree we shall apply B\'{e}zout's
theorem on the projective space ${\bb P}$. By \cite{fu84}
Prop. 8.4, if  $Z\subset {\bb P}$
is a variety and $H$ is a hyperplane such that
$$Z\bigcap H = \bigcup Z_i
$$ where
$Z_i$ are irreducible subvarieties with  $\dim Z_i=\dim Z-1$ then
the degree of $Z$ is given through the formula $\deg Z=\sum
m_i\deg Z_i$ where $m_i$ is the multiplicity of $Z$ and $H$ along
$Z_i$. Furthermore,  by \cite{fu84} Remark 8.2, $m_i=1$
if $Z$  is generically non-singular along $Z_i$ and
generically meets $H$ transversally along $Z_i$.

We first construct a rational map from a
 Schubert variety
of  the Grassmannian
 ${\rm Grass}(m,n(q+1))$ into $Z_\alpha$.
Consider an $m\times n(q+1)$ full rank matrix
$$
Q\in {\rm Grass}(m,n(q+1)).
$$
Let the Pl\"{u}cker coordinate of
$Q\in{\rm Grass}(m,n(q+1))\subset\tilde{\bb P} $ be
$x=(x_i)$ where  $i=(i_1,\dots,i_m)$
and $x_i$ is the $m\times m$ minor of $Q$ consisting of the
$i_1$th through $i_m$th columns. Let
$$
SC_i=\{x\in S_i|x_i\neq 0\}
$$
where $S_i$ is the Schubert variety defined by (\ref{svariety}) with
$n$ replaced by $n(q+1)$.
Then $SC_i$ is a cell and
$$SC_i\cap SC_j=\emptyset,\mbox{ if }\ i\neq j\mbox{ and }
$$
$$
S_i=\bigcup_{j\in I(m), \atop j\leq i} SC_j.
$$

Each $Q\in SC_i$ has a unique echelon form
\begin{equation}\label{echelon}
Q=\left[\begin{array}{ccccccccccccccc}
&&&i_1&&&&i_2&&\cdots&&i_m&&&\\
\ast &\cdots&\ast&1&0&\cdots&0&0&0&\cdots&0&0&0&\cdots&0\\
\vdots&&\vdots&0&\ast&\cdots&\ast&1&0&\cdots&0&0&0&\cdots&0\\
\vdots&&\vdots&\vdots&\vdots&&\vdots&0&\ast&\cdots&0&0&0&\cdots&0\\
\vdots&&\vdots&\vdots&\vdots&&\vdots&\vdots&\vdots&&
\vdots&\vdots&\vdots&&\vdots\\
\ast&\cdots&\ast&0&\ast&\cdots&\ast&0 &\ast&\cdots&\ast&1&0&\cdots&0
\end{array}\right]
\end{equation}

For
$$Q=[Q_0|Q_1|\cdots|Q_q]\in {\rm Grass}(m,n(q+1))
$$
where $\{Q_l\}$ are $m\times n$ matrices,
and  $$
i=(i_1,\dots,i_m),\ 1\leq i_1<\cdots<i_m\leq
n,$$ let
$$
z_{(i;0)}+z_{(i;1)}s+\cdots+z_{(i;pq)}s^{pq}$$
be the $m\times m$ minor formed
by the $i_1$th through $i_m$th columns of the polynomial matrix
\begin{equation}\label{poly}
Q(s)=Q_0+Q_1s+\cdots+Q_qs^q.
\end{equation}
Define a rational map
$\tau : {\rm Grass}(m,n(q+1))\rightarrow {\bb P}$ by
$$
\tau (Q)=(z_{(i;d)})_{1\le i_1<\cdots<i_m\le n,\ 0\le d\le q}
$$
for all points $Q\in {\rm Grass}(m,n(q+1))$
for which the maximum degree of the minors of $Q(s)$ is at most $q$
and at least one minor is non-zero
(these points form a locally closed subset).
Since each $z_{(i;d)}$ is a linear combination of $m\times m$
minors of $Q$, which in turn are the coordinates of $\tilde{\bb P}$
there exists a linear subspace $\tilde{E}\subset\tilde{\bb P}$
such that
$\tau$ is the restriction to ${\rm Grass}(m,n(q+1))$ of the linear
projection $\tau :\tilde{\bb P}
-\tilde{E}\rightarrow {\bb P}$.
We have
\begin{equation}\label{subset}
\tau(S_\alpha)\subset Z_\alpha.
\end{equation}

For a fixed $\alpha=(\alpha_1,\dots,\alpha_m)\in I(m,p)$,
let
$$
\tilde{U}_\alpha=
\bigcup_{\beta=(\beta_1,\dots,\beta_m)\in I(m,p),\atop
\beta\leq\alpha\ {\rm and}\ \beta_1>\alpha_m-n} SC_\beta
\subset {\rm Grass}(m,n(q+1))
$$
and
$$
U_\alpha=
\bigcup_{\beta=(\beta_1,\dots,\beta_m)\in I(m,p),\atop
\beta\leq\alpha\ {\rm and}\ \beta_1>\alpha_m-n} C_\beta
\subset A^q_{m,p}.
$$
Then $\tilde{U}_\alpha$ and $U_\alpha$ are open sets of  $S_\alpha$
and $\overline{C}_\alpha$, respectively.

Let $\phi$ be defined by
$$
\phi(Q)=Q(s)
$$
where $Q(s)$ is defined by (\ref{poly}).
Then the following diagram commutes:
$$
\begin{array}{c}
\tilde{U}_\alpha \ \stackrel{\textstyle \phi}{\rightarrow} \ U_\alpha \\
\tau \searrow\ \ \swarrow \pi\\
Z_\alpha
\end{array}.
$$
If $Q$ is in the echelon form of (\ref{echelon}), then $Q(s)=\phi (Q)$
is in the echelon form defined in~\cite[Proposition 3.5]{wa94}.
Furthermore, if $T$ is the elementary unimodular row operation
which add an $s^k$ multiple of the $l$th row of $Q(s)$ to the
$r$th row,  then $\phi^{-1}(T(Q(s)))$ is in
$$SC_{(\beta_1,\dots,\beta_{r-1},\beta_{r+1},
\dots,\beta_l+kn)},
$$ i.e.
$$\phi^{-1}(T(Q(s))\not\in \tilde{U}_\alpha.
$$
Therefore $\phi:\tilde{U}_{\alpha} \rightarrow U_\alpha$
is one to one and onto.

\begin{pr}\label{pro41}
$S_\alpha$ and $Z_\alpha$ are birationally equivalent under $\tau$.
\end{pr}

\begin{proof}
Let the open set $U$ be defined by
\begin{equation}\label{open}
U=\{Q\in \tilde{U}_\alpha|
\mbox{the $m\times m$ minors of $\phi(Q)$ are relative prime}\}.
\end{equation}
Since
$\pi=\tau\circ \phi^{-1}$
and $\pi:\phi(U)\rightarrow Z_\alpha$ is one to one~\cite{fo75},
$\tau: U\rightarrow Z_\alpha$ is one to one, which means that
$S_\alpha$ is birationally equivalent to
$Z_\alpha$ (see~\cite{ha77}, Chapter I, Corollary 4.5).
\end{proof}

The following proposition generalizes the classical
Pieri  formula and it is one of the main results of this paper.

\begin{pr} \label{zint}
Let $H_\alpha$ be the hyperplane of ${\bb P}$ defined by setting
$z_\alpha=0$. Then
\begin{equation}\label{intersection}
Z_{\alpha}\bigcap H_{\alpha}
=\bigcup_{\beta\in I(m,p),\atop
\beta<\alpha,\ |\beta|=|\alpha|-1}Z_{\beta}
\end{equation}
and the multiplicity of $Z_\alpha$ and $H_\alpha$ along  $Z_\beta$
is one.\end{pr}

\begin{proof}
(\ref{intersection}) follows from Proposition~\ref{o}. So
the only thing we need to prove is that the multiplicity of
the intersection along each $Z_{\beta}$ is one.
Let
$$\tilde{H}_\alpha=\overline{\tau^{-1}(H_\alpha)}=
\tau^{-1}(H_\alpha)\cup \tilde{E}.$$
Then
$$\tau(\tilde{H}_\alpha)=H_\alpha.
$$
$\tilde{H}_\alpha$ is a hyperplane in $\tilde{\bb P}$ defined by
$$0=x_\alpha+\mbox{a linear combination of $x_i$'s
with $i\not\leq \alpha$}.
$$
So
$$S_\alpha\cap \tilde{H}_\alpha=\{x\in S_\alpha|x_\alpha=0\}
=\bigcup_{i\in I(m),\atop
i<\alpha,\ |i|=|\alpha|-1} S_i.
$$

$\tau$
restricted to the open set $U$ defined by~(\ref{open}) is an
isomorphism into an open subset of $Z_\alpha$ and
$U\bigcap S_\beta\neq \emptyset$. So $\tau$ is
a birational isomorphism between $S_\alpha$ and $Z_\alpha$ and
between $S_\beta$ and $Z_\beta$ respectively.
Now, by Pieri's formula~\cite{kl76} applied on the
Grassmannian, the intersection multiplicity of $S_\alpha$ and
$\tilde{H}_\alpha$ along $S_\beta$ is one. By \cite{fu84}~Remark~8.2
and~\ref{pro41},
 the intersection multiplicity of $Z_\alpha$ and
 $H_ \alpha$ along
$Z_\beta$ is also one.
\end{proof}

By B\'{e}zout's Theorem (see discussion following Prop.~8.4~\cite{fu84})
we have the following:
\begin{lem}           \label{deglem}
$$\deg Z_\alpha=
\sum_{\beta\in I(m,p),\atop
\beta<\alpha,\ |\beta|=|\alpha|-1} \deg Z_\beta.
$$ \end{lem}

\begin{th} \label{mainth1}
The degree of $Z_\alpha$ is equal to the number of maximal totally
ordered subsets
of
\begin{equation}\label{indexset}
I_\alpha=\{\beta\in I(m,p)|\beta\leq \alpha\}.
\end{equation}
\end{th}

\begin{proof}
We use induction on the dimension $|\alpha|$ of the
variety $Z_\alpha$. When $|\alpha|=0$, $Z_\alpha=Z_{(1,\dots,p)}$
is a point and its degree is $1$. Assume that the degree of
$Z_\beta$ is  equal to the number of maximal totally
ordered subsets of $I_\beta$
for any $|\beta|=|\alpha|-1$.
Notice that a set $I$ is a maximal totally ordered  subset of
$I_\beta$ for some $\beta$ covered by $\alpha$ if and only if
$I\cup\{\alpha\}$ is a maximal totally ordered subset of $I_\alpha$.
By Lemma~\ref{deglem},
\begin{eqnarray*}
\deg Z_\alpha&=&
\displaystyle{\sum_{\beta\in I(m,p),\atop
\beta< \alpha,\ |\beta|=|\alpha|-1}}
\deg Z_\beta \\
&=&\displaystyle{\sum_{\beta\in I(m,p),\atop
\beta< \alpha,\ |\beta|=|\alpha|-1}}
\mbox{\# of totally ordered subsets of $I_\beta$}\\
&=&\mbox{\# of totally ordered subsets of $I_\alpha$}.
\end{eqnarray*}
\end{proof}

\begin{co}
The degree of $K^q_{m,p}$ is equal to the number of maximal totally
ordered subsets of
$$I_{(p+1,p+2,\dots,n;q)}.$$
\end{co}

In the next section we shall show that the number of maximal totally ordered
subsets is also measured by the formula conjectured in Physics. In order
to show this equivalence, we find it useful to give an abstract
characterization of a function $d(\alpha_1,\dots,\alpha_m)$ that measures the
degree of $Z_\alpha$ as follows:

\begin{co}
The function
\begin{equation}
d(\alpha_1,\dots,\alpha_m):=\deg Z_\alpha
\end{equation}
is the unique solution of the partial recurrence relation
\begin{equation}                 \label{recurr}
d(\alpha_1,\dots,\alpha_m)=\sum_{l=1}^m
d(\alpha_1,\dots,\alpha_l-1,\dots,\alpha_m)
\end{equation}
 subject to the boundary conditions
\begin{eqnarray}
d(\dots,k,k,\dots)&=&0,\label{bound1}\\
d(k,\dots,k+n)&=&0.\label{bound2}
\end{eqnarray}
and subject to the initial conditions
\begin{eqnarray}
d(1,2,\dots,m)&=&1,\label{ini1}\\
d(0,\dots,\alpha_m)&=&0\hspace{2em}
\mbox{ for $\alpha_m<n$.} \label{ini2}
\end{eqnarray}
\end{co}

In~\cite{ra94b} we use  Theorem~\ref{mainth1} and Remark~\ref{remarki}
to give explicit formulas for these degrees.

\section{The conjecture from Physics}
\setcounter{equation}{0}
In this section we shall first give a precise formulation of the
conjecture.
We then interpret the degrees of the subvarieties
$Z_i^d$ that we have just computed through this conjecture. Finally,
we prove that our formula for their degrees
agrees with the conjecture.  We shall not address the physics
behind the conjecture at
all. We refer the reader to~\cite{wi94}~\cite{va92} for a
discussion of the physical aspects.

The cohomology ring of the Grassmannian is generated by the Chern classes
$X_1,\dots ,X_m$ of the canonical subbundle $S$ on ${\rm Grass}(m,n)$.
The complex codimension of the class $X_i$ is $i$ and the
cohomology class $X_m^p$ is Poincar\`e dual to the class of a point.
Let us fix a point $x\in\PP^1$. Then there is an evaluation
map $\psi : X'\to
{\rm Grass}(m,n)$ that sends a map to its value at $x$.
Now any cohomology class on ${\rm Grass}(m,n)$ can be pulled back via
$\psi$ to  $X$. The conjecture predicts
the intersection products of such classes. The first problem with this
is that $X'$ is not a compact space, so one has to interpret these
products on the compactification. But the introduction of the compactification
introduces other questions, namely, do these intersection products
depend on the compactification? In other words, does the boudary component
added to $X'$ change the intersection product.These
questions are somewhat subtle in the case of maps from Riemann
surfaces of higher genus to the Grassmanian
and have been dealt with in~\cite{be93}. In our particular case of maps
from $\PP^1$ to ${\rm Grass}(m,n)$, the compactification chosen by us, namely
the Quot scheme $X$, is a smooth, irreducible variety of
dimension $mp+nq$ for all $m,\, n$ and $q$~\cite{st87}. Further, there is
a universal bundle $\tilde S$ over the Quot scheme $X$ that
extends the pullback $\psi^*(S)$ on $X'$ to all of $X$.
Thus, as in~(\cite{be93}, Section 5.1),
for any set of integers $a_i$ such that $ia_i=mp+nq$
 we can define the
intersection products of $\psi^*(X_i)$ unambiguously through
the following definition:
\begin{equation}\label{intdef}
<\psi^*X_1^{a_1}\dots \psi^*X_m^{a_m}>:=
<c_1^{a_1}\cdots c_m^{a_m}>\end{equation}
where $c_i$ is the $i^{\rm th}$ Chern class of the bundle
$\tilde S$ on $X$.

Let $q_1\dots ,q_m$ be the Chern roots of the canonical subbundle $S$ above,
so $X_i$ is the $i^{\rm th}$ elementary symmetric function
in the $q_j$. Let
$$W=\sum_{i=1}^m(\frac{q_i^{n+1}}{n+1}+(-1)^m q_i).$$
This is called the Landau-Ginzburg potential by physicists. Since
$W$ is a symmetric polynomial in the $q_i$,  it can
be expressed as a polynomial $W=W(X_1,\dots ,X_m)$
in the elementary symmetric functions.
Let
\begin{equation}
h(X_1,\dots ,X_m):=\det[\frac{\partial^2 W}{\partial X_i \partial X_j}]
\end{equation}
be the determinant
of the Hessian of $W(X_1,\dots ,X_m)$. If $a$ is a
multiindex as above
then the intersection product
$<\psi^*X_1^{a_1}\cdots \psi^*X_m^{a_m}>$ is an integer and the conjecture
says that this number is computed by the formula
\begin{equation}\label{conj}
<\psi^*X_1^{a_1}\cdots \psi^*X_m^{a_m}>=(-1)^{m(m-1)/2}
\sum_{dW=0} \frac{X_1^{a_1}\dots X_m^{a_m}}{h}.\end{equation}
The summation in the above formula is over the finite number
of critical points of the function $W(X_1,\dots ,X_m)$.

We will now interpret the degrees that we have computed in terms
of these intersection products. Recall
that $\kappa$ is the map from $X$ to $\PP$ defined by the
line bundle $\wedge^m (\tilde S)$ on $X$ whose Chern class
is $c_1(\tilde S)$. Also
$\kappa(\overline{\psi^{-1}(S_i)})=Z^q_i$. Thus the degree of the
subvariety $Z^q_i\subset \PP$ is given by
$<\psi^*(s_i)\cdot (\psi^*(X_1))^
{{\rm dim}Z^q_i}>$, where
 $s_i(X_1,\dots ,X_m)$
is the Schubert cocycle Poincar\`e dual to $S_i$.
 Further, for $d\le q$, by a
similar argument, the degree of $Z^d_i=<\psi^*(s_i)\cdot (\psi^*(X_1))^
{{\rm dim}Z^d_i}>$.

\begin{th}\label{phythm}
The degree of the subvariety $Z^d_i$ as given by Theorem~\ref{mainth1} also
equals:
\begin{equation}\label{voila}
\begin{array}{ccc}
\deg Z^d_i&=&<\psi^*(s_i)\cdot (\psi^*(X_1))^
{{\rm dim}Z^d_i}>\\
&=&(-1)^{m(m-1)/2}
\sum_{dW=0} \frac{s_i}{h}{X_1}^{\dim Z^d_i}\end{array}\end{equation}
where $h$ is the Hessian of the Landau-Ginzburg potential and
and $s_i$ is the Schubert cocycle of $S_i\subset{\rm Grass}(m,n)$.
In particular one has
\begin{equation}
\deg K^q_{m,p}=(-1)^{m(m-1)/2}
\sum_{dW=0} \frac{X_1^{mp+nq}}{h}
\end{equation}
\end{th}
In order to establish a proof we will first make
some simplifications in the formulas and we will
reformulate the theorem into an equivalent theorem
dealing with properties of the ring of symmetric
functions ${\cal Z}[q_{1},\ldots ,q_{m}]^{{\cal S}_{m}}$ only.

 Some of these simplifications can also be found
in~\cite{wi94}, page 42.
Observe that the Jacobian
$$\det \left[ {\frac{\partial X_i}{\partial q_j} }\right]$$
for the change of variables from $q_i$ to
$X_j$ is given by the Vandermonde determinant
$\Delta =\prod_{j<k}(q_j-q_k)$
and on the critical points the transformation of the Hessian is
given through
$$\det \left.\left[ {\frac{\partial^2 W}{\partial q_i\partial q_j}}\right]
\right|_{dW=0}
=\det \left[ \frac{\partial^2 W}{\partial X_i\partial X_j}\right]
{\left( {\det \left[ {\frac{\partial X_i}{\partial q_j} }\right] }\right)}^2.$$
Thus the polynomial $h$ in terms of the Chern roots $q_j$ is given by
$$h(q_1,\ldots ,q_m)=\frac{n^m\cdot(q_1\dots q_m)^{n-1}}{\Delta^2}.$$
The Schubert cocycle $s_i$ can be identified with
a Schur symmetric function. For this consider the partition
$$\mu := (p+1-i_1,p+2-i_2,\ldots ,p+m-i_m).$$
Let
\begin{equation}
|\mu|:=\mu_1+\cdots+\mu_m=mp-|i|.
\end{equation}
It is well known that the Schubert cocycle $s_i$ can be identified
with the Schur symmetric function $s_\mu$ and $s_\mu$
has a classical representation due to Jacobi ($\sim$1835)
as a quotient of two alternating functions
resulting in a symmetric function:
\begin{equation}        \label{jacobi}
 s_{\mu}(q_1,\dots ,q_m)  \;  =  \;
\frac{\det [q_{i}^{\mu_{j}+m-j}]}{\det[q_{i}^{m-j}]},
i,j=1,\ldots ,m.
\end{equation}
Also, $\Delta=0$ if $q_i=q_j$ so the summation in~(\ref{voila}) is over
all subsets $I$ consisting of $m$ distinct roots of the
the polynomial
$z^n+(-1)^{m}$. On $I$: $(q_1\cdots q_m)^n=1$, so
\begin{equation}
h(q_1,\ldots ,q_m)=\frac{n^m}{(q_1\cdots q_m)\Delta^2}
\end{equation}
Thus to prove Theorem~\ref{phythm} it suffices to
prove the following equivalent theorem:
\begin{th}
The degree of the subvariety $Z^d_i$ as given by Theorem~\ref{mainth1} also
equals:
\begin{equation}       \label{degform}
\deg Z^d_i=\frac{(-1)^{m(m-1)/2}}{n^m}
\sum_{I} (q_1\cdots q_m)
(q_1+\cdots+q_m)^{mp-|\mu |+nq}
\Delta^2 s_\mu
\end{equation}
In particular one has
$$\deg K^q_{m,p}=\frac{(-1)^{m(m+1)/2}}{n^m}
\sum_{I} (q_1\cdots q_m)\left( {\prod_{i<j}(q_i-q_j)^2}\right)
{(q_1+\cdots+q_m)}^{mp+nq}.$$
\end{th}
The proof we will present is purely combinatorial,
 and is just based
on  Proposition~\ref{zint}. In order to establish
the proof we will show that\r{degform} satisfies
the recurrence relation\r{recurr}, the boundary
conditions\r{bound1} and\r{bound2} and
the initial  conditions\r{ini1},\r{ini2}.

\begin{proof}
{}From Jacobi's identity\r{jacobi} and the fact that
$q_i^{n}=(-1)^{m+1}$ it is clear that
formula\r{degform} satisfies both boundary
conditions\r{bound1} and\r{bound2}.
Note that if we let
$$D(\alpha_1,\dots,\alpha_m):=\det[q_i^{n-\alpha_j}],$$
then
$$\det[q_i^{\mu_j+m-j}]=D(i_1,\dots,i_m)=D(\alpha_1,\dots,\alpha_m)$$
where $\alpha=(i;d)$ is defined by~(\ref{ialpha}). So the
recurrence property\r{recurr}
follows from the Pieri type formula:
$$(q_1 +\dots +q_m)\cdot D(\alpha_1,\dots,\alpha_m)=\sum_j
D(\alpha_1,\dots,\alpha_j-1,\dots,\alpha_m).$$

In order to complete the proof it is therefore enough
to show that the initial conditions\r{ini1} and\r{ini2}
are satisfied. Equivalently we have to show that
$$\frac{(-1)^{m(m-1)/2}}{n^m}
\sum_{I}
(q_1\ldots q_m) \Delta^2 s_\mu=\left\{\begin{array}{ll}
1&\mbox{if $\mu=(p^m)$}\\
0&\mbox{if $\mu_1> p$ and $|\mu|=mp$}\end{array}\right. .
$$

We will treat both cases simultaneously.
Note that
$$\Delta=\prod_{j<k} (q_j-q_k)=\det
\left[ \begin{array}{ccc}
1&\ldots &1 \\
\vdots & & \vdots \\
{q_m}^{m-1}&\ldots & {q_1}^{m-1}
\end{array} \right]$$
and
$$
s_{\mu}(q_1,\dots ,q_m)=
\frac{\det [q_{i}^{\mu_{j}+m-j}]}{\Delta}.$$
So we have
\begin{eqnarray}
&&
\frac{(-1)^{m(m-1)/2}}{n^m}
\sum_{I}
(q_1\cdots q_m)\Delta^2 s_\nu\nonumber \\
&=&
\frac{(-1)^{m(m-1)/2}}{n^m}
\sum_{I}
\Delta \det [q_i^{\mu_j+m+1-j}]\nonumber \\
&=&                            \label{product}
\frac{1}{n^m}
\sum_{I}\det
\left[ \begin{array}{ccc}
1&\ldots &1 \\
\vdots & & \vdots \\
{q_1}^{m-1}&\ldots & {q_m}^{m-1}
\end{array} \right]
\left[ \begin{array}{ccc}
{q_1}^{\mu_1+m}&\ldots &{q_1}^{\mu_m+1} \\
\vdots & & \vdots \\
{q_m}^{\mu_1+m}&\ldots & {q_m}^{\mu_m+1}
\end{array} \right]
\end{eqnarray}
Since the summation over the index set $I$ involves
all roots of the polynomial $z^n+(-1)^m$ we can view
the right hand side of the last expression  as
 a symmetric polynomial in
${\cal Z}[y_{1},\ldots ,y_n]^{{\cal S}_n}$,
where $\{ y_{1},\ldots ,y_n\} $ represent all roots
of $z^n+(-1)^m$. This fact is most conveniently
expressed by the Cauchy Binet formula, i.e. the
expression in\r{product} is equal to
\begin{eqnarray}
&=&
\frac{1}{n^m}
\det
\left[ \begin{array}{ccc}
1&\ldots &1 \\
\vdots & & \vdots \\
{y_1}^{m-1}&\ldots & {y_n}^{m-1}
\end{array} \right]
\left[ \begin{array}{ccc}
{y_1}^{\mu_1+m}&\ldots &{y_1}^{\mu_m+1} \\
\vdots & & \vdots \\
{y_n}^{\mu_1+m}&\ldots & {y_{n}}^{\mu_m+1}
\end{array} \right]  \nonumber  \\
&=&
\frac{1}{n^m}
\det
\left[ \begin{array}{cccc}
p(\mu_1 +m)&p(\mu_2 +m-1)&\ldots &p(\mu_m+1) \\
p(\mu_1 +m+1)&p(\mu_2 +m)&\ldots &p(\mu_m+2) \\
\vdots &  & & \vdots \\
p(\mu_1 +2m-1)&p(\mu_2 +2m-2)&\ldots &p(\mu_m+m)
\end{array} \right]
\end{eqnarray}
where $p(k):=\sum_i y_i^k$ is the $k$th power symmetric
function in
${\cal Z}[y_{1},\ldots ,y_{n}]^{{\cal S}_n}$.

Note that
$$
p(k)=\left\{\begin{array}{ll}
n(-1)^{q(m-1)}&\mbox{if $k=qn$}\\
0&\mbox{otherwise}\end{array}\right.$$
Assume now that $\mu =(p^m)$. Then one verifies that
the last  expression evaluates  to
a diagonal matrix with all diagonal entries equal to
$(n)(-1)^{m-1}$.
But this just means that\r{product} correctly evaluates to 1.
 On the other hand
if $\mu \neq (p^m)$ and $|\mu | =mp$ we have  $\mu_m < p$.
But then the last column in the matrix
$\left[ { p(\mu_j+m+i-j)}\right] $ is zero what completes the proof.
\end{proof}

\ifx\undefined\bysame
\newcommand{\bysame}{\leavevmode\hbox to3em{\hrulefill}\,}
\fi
 \end{document}